\magnification 1200 \input miniltx \input graphicx.sty \input amssym \input pictex \font \bbfive = bbm5 \font \bbeight =
bbm8 \font \bbten = bbm10     \font \eightbf = cmbx8
\font \eighti = cmmi8 \skewchar \eighti = '177 \font \eightit = cmti8 \font \eightrm = cmr8 \font \eightsl = cmsl8 \font
\eightsy = cmsy8 \skewchar \eightsy = '60 \font \eighttt = cmtt8 \hyphenchar \eighttt = -1   \font \sixi = cmmi6 \skewchar \sixi = '177 \font \sixrm = cmr6 \font \sixsy = cmsy6 \skewchar \sixsy =
'60 \font \tensc = cmcsc10   \scriptfont \bffam = \bbeight \scriptscriptfont
\bffam = \bbfive \textfont \bffam = \bbten \newskip \ttglue \def \eightpoint {\def \rm {\fam 0 \eightrm }\relax
\textfont 0= \eightrm \scriptfont 0 = \sixrm \scriptscriptfont 0 = \fiverm \textfont 1 = \eighti \scriptfont 1 = \sixi
\scriptscriptfont 1 = \fivei \textfont 2 = \eightsy \scriptfont 2 = \sixsy \scriptscriptfont 2 = \fivesy \textfont 3 =
\tenex \scriptfont 3 = \tenex \scriptscriptfont 3 = \tenex \def \it {\fam \itfam \eightit }\relax \textfont \itfam =
\eightit \def \sl {\fam \slfam \eightsl }\relax \textfont \slfam = \eightsl \def \bf {\fam \bffam \eightbf }\relax
\textfont \bffam = \bbeight \scriptfont \bffam = \bbfive \scriptscriptfont \bffam = \bbfive \def \tt {\fam \ttfam
\eighttt }\relax \textfont \ttfam = \eighttt \tt \ttglue = .5em plus.25em minus.15em \normalbaselineskip = 9pt \def \MF
{{\manual opqr}\-{\manual stuq}}\relax \let \sc = \sixrm \let \big = \eightbig \setbox \strutbox = \hbox {\vrule
height7pt depth2pt width0pt}\relax \normalbaselines \rm } \def \setfont #1{\font \auxfont =#1 \auxfont } \def \withfont
#1#2{{\setfont {#1}#2}} \def \TRUE {Y} \def \FALSE {N} \def \EMPTY {} \def \ifundef #1{\expandafter \ifx \csname
#1\endcsname \relax } \def \undefrule {\kern 2pt \vrule width 5pt height 5pt depth 0pt \kern 2pt} \def \UndefLabels {}
\def \possundef #1{\ifundef {#1}\undefrule {\eighttt #1}\undefrule \global \edef \UndefLabels {\UndefLabels #1\par }
\else \csname #1\endcsname \fi } \newcount \secno \secno = 0 \newcount \stno \stno = 0 \newcount \eqcntr \eqcntr = 0
\ifundef {showlabel} \global \def \showlabel {\FALSE } \fi \ifundef {auxwrite} \global \def \auxwrite {\TRUE } \fi
\ifundef {auxread} \global \def \auxread {\TRUE } \fi \def \define #1#2{\global \expandafter \edef \csname #1\endcsname
{#2}} \long \def \error #1{\medskip \noindent {\bf ******* #1}} \def \fatal #1{\error {#1\par Exiting...}\end } \def
\advseqnumbering {\global \advance \stno by 1 \global \eqcntr =0} \def \current {\ifnum \secno = 0 \number \stno \else
\number \secno \ifnum \stno = 0 \else .\number \stno \fi \fi } \begingroup \catcode `\@=0 \catcode `\\=11 @global @def
@textbackslash {\} @endgroup \def \space { } \def \rem #1{\vadjust {\vbox to 0pt{\vss \hfill \raise 3.5pt \hbox to 0pt{
#1\hss }}}} \font \tiny = cmr6 scaled 800 \def \deflabel #1#2{\relax \if \TRUE \showlabel \rem {\tiny [#1]}\fi \ifundef
{#1PrimarilyDefined}\relax \define {#1}{#2}\relax \define {#1PrimarilyDefined}{#2}\relax \if \TRUE \auxwrite \immediate
\write 1 {\textbackslash newlabel {#1}{#2}}\fi \else \edef \old {\csname #1\endcsname }\relax \edef \new {#2}\relax \if
\old \new \else \fatal {Duplicate definition for label ``{\tt #1}'', already defined as ``{\tt \old }''.}\fi \fi }  \def \label #1 {\deflabel {#1}{\current }} \def \equationmark #1 {\ifundef
{InsideBlock} \advseqnumbering \eqno {(\current )} \deflabel {#1}{\current } \else \global \advance \eqcntr by 1 \edef
\subeqmarkaux {\current .\number \eqcntr } \eqno {(\subeqmarkaux )} \deflabel {#1}{\subeqmarkaux } \fi } \def \split
#1.#2.#3.#4;{\global \def \parone {#1}\global \def \partwo {#2}\global \def \parthree {#3}\global \def \parfour {#4}}
\def \NA {NA} \def \ref #1{\split #1.NA.NA.NA;(\possundef {\parone }\ifx \partwo \NA \else .\partwo \fi )}  \newcount \bibno \bibno = 0  \def \Bibitem #1 #2; #3;
#4 \par {\smallbreak \global \advance \bibno by 1 \item {[\possundef {#1}]} #2, {``#3''}, #4.\par \ifundef
{#1PrimarilyDefined}\else \fatal {Duplicate definition for bibliography item ``{\tt #1}'', already defined in ``{\tt
[\csname #1\endcsname ]}''.}  \fi \ifundef {#1}\else \edef \prevNum {\csname #1\endcsname } \ifnum \bibno =\prevNum
\else \error {Mismatch bibliography item ``{\tt #1}'', defined earlier (in aux file ?) as ``{\tt \prevNum }'' but should
be ``{\tt \number \bibno }''.  Running again should fix this.}  \fi \fi \define {#1PrimarilyDefined}{#2}\relax \if \TRUE
\auxwrite \immediate \write 1 {\textbackslash newbib {#1}{\number \bibno }}\fi } \def \jrn #1, #2 (#3), #4-#5;{{\sl #1},
{\bf #2} (#3), #4--#5} \def \Article #1 #2; #3; #4 \par {\Bibitem #1 #2; #3; \jrn #4; \par } \def \references
{\begingroup \bigbreak \eightpoint \centerline {\tensc References} \nobreak \medskip \frenchspacing } \catcode `\@=11
\def \c@itrk #1{{\bf \possundef {#1}}} \def \c@ite #1{{\rm [\c@itrk{#1}]}} \def \sc@ite [#1]#2{{\rm [\c@itrk{#2}\hskip
0.7pt:\hskip 2pt #1]}} \def \du@lcite {\if \pe@k [\expandafter \sc@ite \else \expandafter \c@ite \fi } \def \cite
{\futurelet \pe@k \du@lcite } \catcode `\@=12 \def \Headlines #1#2{\nopagenumbers \headline {\ifnum \pageno = 1 \hfil
\else \ifodd \pageno \tensc \hfil \lcase {#1} \hfil \folio \else \tensc \folio \hfil \lcase {#2} \hfil \fi \fi }} \def
\title #1{\medskip \centerline {\withfont {cmbx12}{\ucase {#1}}}}   \long \def \Quote #1\endQuote {\begingroup \leftskip 35pt
\rightskip 35pt \parindent 17pt \eightpoint #1\par \endgroup } \long \def \Abstract #1\endAbstract {\vskip 1cm \Quote
\noindent #1\endQuote }   \def \Note #1{\footnote {}{\eightpoint #1}} \def \Date #1 {\Note {\it Date: #1.}}
\newcount \auxone \newcount \auxtwo \newcount \auxthree \def \currenttime {\auxone =\time \auxtwo =\time \divide \auxone
by 60 \auxthree =\auxone \multiply \auxthree by 60 \advance \auxtwo by -\auxthree \ifnum \auxone <10 0\fi \number
\auxone :\ifnum \auxtwo <10 0\fi \number \auxtwo } \def \today {\ifcase \month \or January\or February\or March\or
April\or May\or June\or July\or August\or September\or October\or November\or December\fi \space \number \day , \number
\year }  \def \hojeExtenso {\number \day \ de \ifcase \month
\or janeiro\or fevereiro\or mar\c co\or abril\or maio\or junho\or julho\or agosto\or setembro\or outubro\or novembro\or
decembro\fi \ de \number \year }  \def \part #1#2{\vfill \eject \null
\vskip 0.3\vsize \withfont {cmbx10 scaled 1440}{\centerline {PART #1} \vskip 1.5cm \centerline {#2}} \vfill \eject }
  \def \fix {\smallskip \noindent $\blacktriangleright $\kern 12pt} \def \iskip {\medskip
\noindent }    \def \ucase #1{\edef \auxvar
{\uppercase {#1}}\auxvar } \def \lcase #1{\edef \auxvar {\lowercase {#1}}\auxvar } \def \emph #1{{\it #1}} \def \section
#1 \par {\global \advance \secno by 1 \stno = 0 \goodbreak \bigbreak \noindent {\bf \number \secno .\enspace #1.}
\nobreak \medskip \noindent } \def \state #1 #2\par {\begingroup \def \InsideBlock {} \medbreak \noindent
\advseqnumbering {\bf \current .\enspace #1.\enspace \sl #2\par }\medbreak \endgroup } \def \definition #1\par {\state
Definition \rm #1\par } \newcount \CloseProofFlag \def \closeProof {\eqno \endproofmarker \global \CloseProofFlag =1}
 \long \def \Proof #1\endProof {\begingroup \def
\InsideBlock {} \global \CloseProofFlag =0 \medbreak \noindent {\it Proof.\enspace }#1 \ifnum \CloseProofFlag =0 \hfill
$\endproofmarker $ \looseness = -1 \fi \medbreak \endgroup } \def \quebra #1{#1 $$$$ #1} \def \explica #1#2{\mathrel
{\buildrel \hbox {\sixrm #1} \over #2}} \def \explain #1#2{\explica {\ref {#1}}{#2}}  \def \=#1{\explain {#1}{=}} \def \pilar #1{\vrule height #1 width 0pt}  \newcount \fnctr \fnctr = 0 \def \fn #1{\global \advance \fnctr by 1 \edef \footnumb
{$^{\number \fnctr }$}\relax \footnote {\footnumb }{\eightpoint #1\par \vskip -10pt}} \def \text #1{\hbox {#1}}    \def \Item #1{\smallskip \item {{\rm #1}}}
\newcount \zitemno \zitemno = 0 \def \izitem {\global \zitemno = 0}  \def \zitemplus {\global
\advance \zitemno by 1 \relax } \def \rzitem {\romannumeral \zitemno } \def \rzitemplus {\zitemplus \rzitem } \def
\zitem {\Item {{\rm (\rzitemplus )}}}   \def \zitemmark #1 {\deflabel
{#1}{\current .\rzitem }\relax {\def \showlabel {\FALSE }\deflabel {Local#1}{\rzitem }}} \def \iItemmark #1 {\zitemmark
{#1} } \def \zitemlbl #1 {\zitem \deflabel {#1}{\current .\rzitem }\relax {\def \showlabel {\FALSE }\deflabel
{Local#1}{\rzitem }}} \newcount \nitemno \nitemno = 0  \def \nitem {\global \advance \nitemno
by 1 \Item {{\rm (\number \nitemno )}}} \newcount \aitemno \aitemno = -1 \def \boxlet #1{\hbox to 6.5pt{\hfill #1\hfill
}} \def \iaitem {\aitemno = -1} \def \aitemconv {\ifcase \aitemno a\or b\or c\or d\or e\or f\or g\or h\or i\or j\or k\or
l\or m\or n\or o\or p\or q\or r\or s\or t\or u\or v\or w\or x\or y\or z\else zzz\fi } \def \aitem {\global \advance
\aitemno by 1\Item {(\boxlet \aitemconv )}} \def \aitemmark #1 {\deflabel {#1}{\aitemconv }}  \def \Case #1:{\medskip \noindent {\tensc Case #1:}} \def \<{\left \langle \vrule width 0pt depth 0pt
height 8pt } \def \>{\right \rangle } \def \({\big (} \def \){\big )}  \def \and {\hbox {\quad
and \quad }} \def \calcat #1{\,{\vrule height8pt depth4pt}_{\,#1}}   \def \IMPLY {\kern 7pt \Rightarrow \kern 7pt} \def \for #1{\quad \forall \,#1} \def
\endproofmarker {\square } \def \"#1{{\it #1}\/}  \def \inv {^{-1}} \def \*{\otimes }
\def \caldef #1{\global \expandafter \edef \csname #1\endcsname {{\cal #1}}} \def \mathcal #1{{\cal #1}} \def \bfdef
#1{\global \expandafter \edef \csname #1\endcsname {{\bf #1}}} \bfdef N \bfdef Z \bfdef C \bfdef R  \def \exists {\mathchar "0239\kern 1pt } \if \TRUE \auxread \IfFileExists {\jobname .aux}{\input \jobname
.aux}{\null } \fi \if \TRUE \auxwrite \immediate \openout 1 \jobname .aux \fi \def \close {\if \EMPTY \UndefLabels \else
\message {*** There were undefined labels ***} \iskip ****************** \ Undefined Labels: \tt \par \UndefLabels \fi
\if \TRUE \auxwrite \closeout 1 \fi \par \vfill \supereject \end }  \newcount \pnumber \def \eqnarray
#1#2{\centerline {\vbox {\halign {$##$&\enspace $#1$\enspace $##$\vrule width 0pt height 10pt depth 5pt\hfil \cr #2}}}}
\def \clauses #1{\def \crr {\vrule width 0pt height 10pt depth 5pt\cr }\left \{ \matrix {#1}\right .}  \def \cl #1 #2 #3
{#1, & \hbox {#2 } #3\hfill \crr } \def \beginmypix x from #1 to #2 y from #3 to #4 xscale #5 yscale #6 { \begingroup
\def \LeftLimit {#1} \def \RightLimit {#2} \def \TopLimit {#3} \def \BotLimit {#4} \noindent \hfill \beginpicture
\setcoordinatesystem units <#5, #6> \setplotarea x from {\LeftLimit } to {\RightLimit }, y from {\TopLimit } to
{\BotLimit } \put {\null } at {\LeftLimit } {\TopLimit } \put {\null } at {\RightLimit } {\BotLimit } } \def \endmypix
{\endpicture \hfill \null \endgroup }  \newcount \ans \def \slide #1 #2 #3
#4 {\relax \ans =#2 \advance \ans by -#1 \multiply \ans by #3 \divide \ans by 100 \advance \ans by #1 \global
\expandafter \edef \csname #4\endcsname {\number \ans }} \newcount \ax \newcount \ay \newcount \bx \newcount \by
\newcount \mx \newcount \my \newcount \aux \newcount \afasta \afasta = 20 \def \flexa #1 #2 #3 #4 {\aux = 100 \advance
\aux by -\afasta \slide #1 #3 {\number \afasta } ax \slide #1 #3 {\number \aux } bx \slide #2 #4 {\number \afasta } ay
\slide #2 #4 {\number \aux } by \arrow <0.15cm> [0.25,0.75] from {\ax } {\ay } to {\bx } {\by } \slide #1 #3 50 mx
\slide #2 #4 50 my } \def \map #1 #2 #3 #4 #5 #6 {\put {#1} at #2 #3 \put {#4} at #5 #6 \flexa #2 #3 #5 #6 } \def \lbl
<#1,#2> #3{\put {#3} <#1,#2> at {\mx } {\my } } \def \IFF {\mathrel {\kern 2pt \Leftrightarrow \kern 2pt}} \def \sub
#1{_{\scriptscriptstyle #1}} \def \1{{\bf 1}} \def \Markov {M_\Lambda } \def \wkautomat {weak $k$-automat} \font \bigrm
= cmbx12 scaled 1440 \font \bigmath = cmmi12 scaled 1440

% Higher rank graphs, k-subshifts and k-automata
% R. Exel and B. Steinberg

  \bigrm
  \centerline {Higher rank graphs, {\bigmath k}-subshifts}
  \medskip
  \centerline {and {\bigmath k}-automata}
  \bigskip
  \tensc
  \centerline {
  R. Exel\footnote {$^{\ast }$}{\eightrm Universidade Federal de Santa Catarina and University of Nebraska -- Lincoln.}
and
  B. Steinberg\footnote {$^{\ast\ast} $}{\eightrm The City University of New York.}}
  \rm

  \Abstract Given a $k$-graph $\Lambda $ we construct a Markov space $\Markov $, and a collection of $k$ pairwise
commuting cellular automata on $\Markov $, providing for a factorization of Markov's shift.  Iterating these maps we
obtain an action of ${\bf N}^k$ on $\Markov $ which is then used to form a semidirect product groupoid $\Markov \rtimes
{\bf N}^k$.  This groupoid turns out to be identical to the path groupoid constructed by Kumjian and Pask, and hence its
C*-algebra is isomorphic to the higher rank graph C*-algebra of $\Lambda $.
  \endAbstract

\section Introduction

Given a row-finite $k$-graph $\Lambda $ (see \cite {KP}), Kumjian and Pask have constructed a \emph {path space}
$\Lambda ^\infty $ and an action of ${\bf N}^k$ on $\Lambda ^\infty $ such that the C*-algebra of the \emph {path
groupoid} $\Lambda ^\infty \rtimes {\bf N}^k$ is canonically isomorphic to the corresponding higher rank graph
C*-algebra $C^*(\Lambda )$.

Recall from \cite {KP} that a path in $\Lambda $ consists of a map
  $
  x:\Omega _k \to \Lambda ,
  $
  where
  $$
  \Omega _k=\big \{(n, m)\in {\bf N}^k\times {\bf N}^k, \ n\leq m\big \},
  $$
  satisfying suitable conditions.  In the paragraph after \cite [Remarks 2.2]{KP}, the authors observe that each path
$x$ in $\Lambda ^\infty $ is uniquely determined by a very small subset of its values, such as, for example, the values
of the form
  $$
  y(n)\,= \,x\big (\,(n, n, \ldots , n)\, ,\,(n{+}1, n{+}1, \ldots , n{+}1)\,\big ),
  $$
  for every $n\in {\bf N}$. Noting that each $y(n)$ above is an element of $\Lambda $ of degree
  $$
  d\big(y(n)\big)=(1,1, \ldots ,1),
  $$
  we consider the subset $\Sigma $ of $\Lambda $ formed by all elements possessing the above degree.  Viewing $\Sigma $
as an alphabet, in the spirit of Symbolic Dynamics, one may easily see that each $y$ given above is in fact an element
of a certain Markov subspace $\Markov \subseteq \Sigma ^{\bf N}$, the correspondence $x\to y$ in fact being a
homeomorphism from $\Lambda ^\infty $ to $\Markov $.

We thus have two homeomorphic spaces, each carrying an action of a different monoid, namely of ${\bf N}^k$ in case of
$\Lambda ^\infty $, while ${\bf N}$ acts on $\Markov $ by means of iterating Markov's shift.  These actions are
compatible in the sense that the above correspondence is covariant for the action of the submonoid
  $$
  \big \{(n,n,\ldots ,n):n\in {\bf N}\big \}\subseteq {\bf N}^k,
  $$
  on $\Lambda ^\infty $, but one may also use the above homeomorphism to extend the action of ${\bf N}$ on $\Markov $ to
an action of the much larger monoid ${\bf N}^k$.  Alternatively, considering the action of the canonical basis vectors
$e_i$ of $\, {\bf N}^k$ on $\Lambda ^\infty $, we may define continuous maps
  $$
  S_i:\Markov \to \Markov ,
  $$
  for $i=1,\ldots ,k$,
  which commute among themselves, giving a factorization of Markov's shift $S$ on $\Markov $, in the sense that
  $$
  S_1S_2\cdots S_k=S.
  $$

The well known Curtis-Hedlund-Lyndon Theorem in fact states that, when the alphabet is finite, every continuous map
commuting with the shift on $\Sigma ^{\bf N}$ is a \emph {cellular automaton}, given by means of a \emph {sliding block
code}.  Regardless of the size of our alphabet, we indeed show that the $S_i$ above are given by sliding block codes
closely linked to the unique factorization property of $\Lambda $.

Motivated by this example, we introduce the notion of a \emph {\wkautomat on} over a given alphabet $\Sigma $, as being
a $k{+}1$-tuple
  $$
  (Y;\ S_1,S_2,\ldots ,S_k),
  $$
  where $Y$ is a classical subshift (i.e., a closed subset of $\Sigma ^{\bf N}$, invariant under the shift $S$), and the
$S_i$ are pairwise commuting continuous maps from $Y$ to $Y$, providing for a factorization of the shift.

Since the path groupoid $\Lambda ^\infty \rtimes {\bf N}^k$ may be built from nothing more than the information
contained in the action of ${\bf N}^k$ on the path space, one sees that this groupoid is identical in all respects to
the groupoid constructed from the associated \wkautomat on, and hence that the higher rank graph C*-algebra may be
constructed solely based on the the latter.

As an auxiliar gadget we also define a notion of a \emph {$k$-subshift}, as being a closed subset of $\Sigma ^{{\bf
N}^k}$, invariant under the natural action of ${\bf N}^k$, and which is isomorphic to its image under the restriction to
the diagonal.  See \ref {DefKSubsh} for the precise definition.  The relevance of this notion resides in the fact that
it has a more geometrical appeal, while being essentially the same thing as a \wkautomat on, as proved in \ref
{SameThing}.

The adjective ``weak'' above is nothing but a disclaimer highlighting the fact that the $S_i$ involved are not actually
supposed to be cellular automata, although they share with the latter the important property of commuting with the
shift.  In \ref{CharacStrongAutomata} we then improve on \ref{SameThing} by precisely characterizing the $k$-subshifts
giving rise to \wkautomat a involving actual cellular automata.

  The first named author was partially supported by CNPq.
  The second named author thanks the Fulbright Commission for supporting his recent visit to Florianopolis during which
part of the research for this paper was conducted.

\section $P$-subshifts

Let $\Sigma $ be a set, henceforth called the \emph {alphabet}, viewed as a topological space with the discrete
topology.

Given a monoid $P$, we will consider the set $\Sigma ^P$ equipped with the product topology.  Although we will not
assume that $\Sigma $ is finite here, we observe that when $\Sigma $ is finite, then $\Sigma ^P$ is a compact space by
Tychonov's Theorem.

For each $p$ in $P$, we will moreover denote by
  $
  \theta _p:\Sigma ^P\to \Sigma ^P,
  $
  the map given by
  $$
  \theta _p(\xi )\calcat t = \xi (pt), \for \xi \in \Sigma ^P, \for t\in P.
  $$
  It is then easy to prove that $\theta $ is a right action of $P$ on $\Sigma ^P$, that is,
  $$
  \theta _p\theta _q = \theta _{qp}, \for p, q \in P,
  $$
  henceforth referred to as the \emph {Bernoulli action}, or the \emph {full $P$-shift}, on the alphabet $\Sigma $.

\fix From now on the alphabet $\Sigma $ and the semigroup $P$ will be considered fixed.

\definition \label DefKSubsh
  A \emph {$P$-subshift}\/ is any closed subset $X\subseteq \Sigma ^P$ which is invariant under the Bernoulli action in
the sense that $\theta _p(X)\subseteq X$, for every $p$ in $P$.

We will next describe an important source of $P$-subshifts given in terms of \emph {forbidden patterns}.

\definition
  By a \emph {pattern} we shall mean a pair $(\pi , D_\pi )$, where $D_\pi $ is a finite subset of $P$, and $\pi :D_\pi
\to \Sigma $ is any function.  We shall frequently refer to $\pi $ as a pattern without mentioning $D_\pi $ explicitly.
Given a pattern $\pi $ and an element $\xi $ in $\Sigma ^P$, we will say that \emph {$\pi $ occurs in $\xi $}, provided
there exists some $p_0$ in $P$ such that $\pi (t) = \xi (p_0t)$, for all $t$ in $D_\pi $.

\state Proposition \label ContraPos
  If the pattern $\pi $ occurs in $\theta _p(\xi )$, for some $\xi $ in $\Sigma ^P$, then $\pi $ occurs in $\xi $.

\Proof By hypothesis there exists $p_0$ in $P$ such that
  $$
  \pi (t) = \theta _p (\xi )\calcat {p_0t} = \xi (pp_0t), \for t\in D_\pi ,
  $$
  whence the conclusion.
  \endProof

Given a pattern $\pi $, it is easy to see that the set of all $\xi $ in $\Sigma ^P$ such that $\pi $ occurs in $\xi $ is
open in $\Sigma ^P$.

\state Proposition
  Given a collection $\Pi $ of patterns, let $X_\Pi $ be the set of all elements $\xi $ in $\Sigma ^P$ such that
\underbar {no pattern} in $\Pi $ occurs in $\xi $.  Then $X_\Pi $ is a $P$-subshift.

\Proof
  $X_\Pi $ is closed by the observation made just before the statement, and it is invariant under the Bernoulli action
by the contrapositive of \ref {ContraPos}.  \endProof

The following is a well known result in the theory of classical subshifts.  It is usually stated for finite alphabets,
but it works just as well for infinite ones.

\state Proposition
  If $X$ is any $P$-subshift then there exists a collection $\Pi $ of patterns such that $X=X_\Pi $.

\Proof
  Let $\Pi $ be the collection of all patterns which do not occur in any $\xi $ in $X$.  It is then obvious that
$X\subseteq X_\Pi $, and we next claim that $X$ is dense in $X_\Pi $.  To see this, choose any $\eta $ in $X_\Pi $, and
let $V$ be any open subset of $\Sigma ^P$ containing $\eta $.  Since $\Sigma ^P$ has the product topology, there exists
a finite set $D\subseteq P$ such that
  $$
  \eta \in W:=\{\zeta \in \Sigma ^P: \zeta (t)=\eta (t), \hbox { for all } t\in D\}\subseteq V.
  $$
  Setting $\pi =\eta |_D$, we have that $\pi $ is a pattern, obviously occuring in $\eta $, whence $\pi $ is certainly
not in $\Pi $.  By definition of $\Pi $, it follows that $\pi $ occurs in some $\xi $ in $X$, so there exists $p_0$ in
$P$ such that
  $$
  \pi (t) = \xi (p_0t) = \theta _{p_0}(\xi )\calcat t,
  $$
  for every $t$ in $D$.  This says that $\theta _{p_0}(\xi )\in W\cap X$, proving the desired density, namely that
$X_\Pi \subseteq \overline {X\pilar {8pt}}$.  Since $X$ is closed by hypothesis, the proof is concluded.  \endProof

If $Q$ is a submonoid of $P$, always assumed to share the neutral element, we may consider the \emph {restriction
mapping} \ $\rho \sub Q$ \ from $\Sigma ^P$ to $\Sigma ^Q$, namely
  $$
  \rho \sub Q(\xi ) = \xi |\sub Q, \for \xi \in \Sigma ^P.
  \equationmark RestrMap
  $$
  Clearly $\rho \sub Q$ is a continuous map.

\state Proposition \label RestrictPShift
  Let $X\subseteq \Sigma ^P$ be a $P$-subshift, and let $Q$ be a submonoid of $P$. Then
  \izitem
  \zitem $\rho \sub Q(X)$ is invariant under the Bernoulli action of\/ $Q$,
  \zitem if\/ $\Sigma $ is finite, then $\rho \sub Q(X)$ is a $Q$-subshift.

\Proof
  Letting $\theta '$ denote the full $Q$-shift, observe that, for all $q$ in $Q$, and all $\xi $ in $X$, one has that
  $$
  \theta '_q(\xi |\sub Q) = \theta _q (\xi )|\sub Q,
  $$
  from where it easily follows that $\rho \sub Q(X)$ is invariant under the Bernoulli action of $Q$.  Assuming that
$\Sigma $ is finite, we have that $\Sigma ^P$ is compact, whence so is $X$.  Observing that $\rho \sub Q$ is continuous,
we see that $\rho \sub Q(X)$ is compact, hence closed in $\Sigma ^Q$.  This concludes the proof.  \endProof

We will soon discuss an important class of examples in which $\rho \sub Q(X)$ is closed, even though the alphabet might
be infinite.  It will then follow from \ref {RestrictPShift.i} that $\rho \sub Q(X)$ is a $Q$-subshift.  Incidentally,
we do not have any example in which $\rho \sub Q(X)$ fails to be closed.

\section Cellular Automata

\label CellAutomataSection As before we let $\Sigma $ be any set, which we view as a discrete topological space.  From
now on we shall be concerned with metric aspects, most notably with the notion of uniform continuity, so we shall equip
$\Sigma $ with the metric defined by
  $$
  d(a, b) = \clauses{
  \cl 0 if a=b,
  \cl 1 otherwise, {} }
  $$
  for all $a$ and $b$ in $\Sigma $.
  The topology on $\Sigma $ induced by this metric is clearly the discrete topology.
  For this reason $d$ is sometimes called the \emph{discrete metric}.  However, while many other metrics on $\Sigma $
also induce the discrete topology, we observe that $d$ induces the discrete topology in a \emph{uniform} way, meaning
that there exists $r>0$, namely $r=1/2$ such that for every $a$ in $\Sigma $, the ball centered at $a$ with radius $r$
coincides with the singleton $\{a\}$.  The fact that $r$ does not depend on $a$ is what makes $d$ a \emph{uniformly
discrete} metric.

In this section we shall be concerned with the monoid ${\bf N}$, formed by all natural numbers, including zero, and
hence we will be working with ${\bf N}$-subshifts, also known simply as subshifts.

The most popular metric considered on the product space $\Sigma ^{\bf N}$ (as usual also denoted by $d$, by abuse of
language) is as follows: given $x$ and $y$ in $\Sigma ^{\bf N}$, one puts $d(x,y) = 2^{-k}$, where $k$ is the largest
integer such that $x_i=y_i$, for all $i\leq k$.  If $x=y$, then obviously no such $k$ exists, in which case we set
$d(x,y)=0$.  It is well known that $d$ defines a metric on $\Sigma ^{\bf N}$, which is compatible with the product
topology.

The role of uniform continuity is evidenced by our next result.

\state Lemma \label UnifCont
  \izitem
  \zitem Each projection $p_k:\Sigma ^{\bf N}\to \Sigma $ is uniformly continuous.
  \zitem For every nonempty $X\subseteq \Sigma ^{\bf N}$, and for every uniformly continuous map $\varphi :X\to \Sigma
$, one has that $\varphi $ depends only on finitely many coordinates, meaning that there exists some $k\in {\bf N}$, and
a map $\psi :\Sigma ^{k+1}\to \Sigma $ such that
  $$
  \varphi (x) = \psi (x_0,x_1,\ldots ,x_k), \for x\in X.
  $$

\Proof Regarding the projection $p_k$, and given $\varepsilon >0$, choose $\delta =2^{-k}$.  Then, for every $x$ and $y$
in $\Sigma ^{\bf N}$, with $d(x,y)\leq \delta $, we necessarily have that $x_i=y_i$, for all $i\leq k$, hence
  $$
  d\big(p_k(x),p_k(y)\big) = d(x_k, y_k) = 0 < \varepsilon ,
  $$
  proving (i).  With respect to (ii), let $\delta >0$ be such that
  $$
  d(x,y)<\delta \IMPLY d\big(\varphi (x),\varphi (y)\big)<1/2,
  $$
  for every $x,y\in \Sigma ^{\bf N}$,
  and choose an integer $k$ such that $2^{-k}<\delta $.  Given $x$ and $y$ in $\Sigma ^{\bf N}$ such that
  $$
  (x_0,x_1,\ldots ,x_k) = (y_0,y_1,\ldots ,y_k),
  $$
  we then have that
  $
  d(x,y)\leq 2^{-k}<\delta ,
  $
  so $d\big(\varphi (x),\varphi (y)\big)<1/2$, which clearly implies that $\varphi (x)=\varphi (y)$, since $\Sigma $ has
the 0-1 metric.  This proves that $\varphi (x)$ depends only on $(x_0,x_1,\ldots ,x_k)$.  \endProof

Observe that when $\Sigma $ is finite and $X$ is closed in $\Sigma ^{\bf N}$, then $X$ is compact by \hbox{Tychonov's}
Theorem, so every continuous function on $X$ is necessarily uniformly continuous.  Consequently the conclusion of
\ref{UnifCont.ii} holds for every continuous function $\varphi $.

\bigskip As usual, we denote by $S$ the \emph{shift} on $\Sigma ^{\bf N}$, defined by
  $$
  S(x_0,x_1,x_2, \ldots ) = (x_1,x_2,x_3, \ldots ), \for x=(x_0,x_1,x_2, \ldots ) \in \Sigma ^{\bf N},
  $$
  so that a closed subspace $X\subseteq \Sigma ^{\bf N}$ is an ${\bf N}$-subshift if and only if $X$ is invariant under
$S$ in the sense that $S(X)\subseteq X$.

\def\L{{\cal L}}

\definition Let $X\subseteq \Sigma ^{\bf N}$ be a subshift.
  \izitem
  \zitem The \emph{language} of $X$, denoted $\L(X)$, is the set of all finite words occurring as a contiguous block of
characters in some $x$ in $X$.
  \zitem Given $k\in {\bf N}$, the subset formed by all words in $\L(X)$ of length $k+1$ will be denoted by $\L_k(X)$.
  \zitem By a \emph{(sliding block) code} for $X$ we mean any function $\psi :\L_k(X)\to \Sigma $.  The integer $k$ is
called the \emph{anticipation} of $\psi $.
  \zitem Given a code $\psi $, we define
  $$
  T_\psi :(x_0, x_1, \ldots )\in X\mapsto (y_0, y_1, \ldots )\in \Sigma ^{\bf N},
  $$
  where $y_n= \psi (x_n,x_{n+1},\ldots,x_{n+k})$,
  and $k$ is the anticipation of $\psi $.  One says that $T_\psi $ is the \emph{cellular automaton} associated to the
code $\psi $.

Another reason for our interest in uniform continuity is in order:

\state Proposition \label CellUnif
  Let $X\subseteq \Sigma ^{\bf N}$ be a subshift and let $\psi :\L_k(X)\to \Sigma $ be a code for $X$.  The the cellular
automaton $T_\psi $ is a uniformly continuous map.

\Proof
  Given $\varepsilon >0$, choose $p$ such that $2^{-p}<\varepsilon $.  For $x$ and $x'$ in $X$, let
  $$
  y=T_\psi (x), \and y'=T_\psi (x').
  $$

  If $d(x,x') \leq 2^{-(p+k)}$, we have that $x_i=x'_i$, for all $i\leq p+k$, from where one deduces that $y_i=y'_i$,
for all $i\leq p$, whence
  $$
  d(y,y') \leq 2^{-p} < \varepsilon .
  $$
  This proves that $T_\psi $ is uniformly continuous.
  \endProof

If $\psi $ is defined on $\Sigma ^2$ by $\psi (a,b)=b$, then $T_\psi $ is clearly the shift itself.  From the above we
then deduce the elementary fact that the shift is uniformly continuous.

Given a subshift $X$ and a code $\psi $ for $X$, it easy to prove that
  $$
  T_\psi \big(S(x)\big)=S\big(T_\psi (x)\big), \for x \in X.
  \equationmark Commutativity
  $$

The next result, known as the Curtis-Hedlund-Lyndon Theorem, says that the only uniformly continuous maps commuting with
the shift are the cellular automata.  The corresponding result for the bilateral shift is discussed in
\cite[13.9]{LindMarcus}.

\state Proposition \label Hedlund
  Let $X$ be a subshift and let $T:X\to X$ be a continuous mapping which commutes with $S$.
  Then $T$ is uniformly continuous if and only if $T$ is a cellular automaton.

\Proof The ``if'' part having already been dealt with in \ref{CellUnif}, we move on to the ``only if'' part.  Consider
the function $\varphi :X\to \Sigma $ given by $\varphi =p_0\circ T$, where $p_0$ is the projection on the leftmost
coordinate.  By \ref{UnifCont.i} we have that $\varphi $ is uniformly continuous, and hence
  $$
  p_0\big(T(x)\big) = \psi (x_0,x_1,\ldots , x_k), \for x\in X,
  $$
  for some $k$ and some $\psi :\Sigma ^{k+1}\to \Sigma $, by \ref{UnifCont.ii}.  For every $x$ in $X$ it then follows
that
  $$
  p_n\big(T(x)\big) =
  p_0\big(S^n(T(x))\big) =
  p_0\big(T(S^n(x))\big) \quebra=
  p_0\big(T(x_n, x_{n+1}, \ldots )\big) =
  \psi (x_n,x_{n+1},\ldots,x_{n+k}).
  $$

Noticing that $(x_n,x_{n+1},\ldots,x_{n+k})$ necessarily lies in $\L_k(X)$, and upon restricting $\psi $ to $\L_k(X)$,
we then have that $T$ is the cellular automaton associated to $\psi $.  \endProof

\section Higher rank graphs

Given any integer $k\geq 1$, let $\Lambda $ be a $k$-graph (see \cite {KP}).  That is, $\Lambda$ is a (small) category
equipped with a degree functor $d\colon \Lambda\to {\bf N}^k$, satisfying the \emph{unique factorization property},
namely, if $d(\lambda)=m+n$, then there are
unique $\alpha,\beta\in \Lambda$ with $\lambda=\alpha\beta$ and $d(\alpha)=m$ and $d(\beta)=n$.

  Recall that $\Omega _k$ is the $k$-graph consisting of all pairs $(m,n)\in {\bf N}^k\times {\bf N}^k$ such that $m\leq
n$, equipped with the degree map defined by $d(m,n)=n-m$ and with allowed products $(n,r)(m,n)=(m,r)$.  By definition, a
path in $\Lambda $ is a functor from $\Omega _k$ to $\Lambda $, compatible with the degree maps.  The set of all paths
in $\Lambda $ is denoted by $\Lambda ^\infty $.

Recall from \cite [Definitions 2.1]{KP} that, for each $p$ in ${\bf N}^k$, one defines $\sigma ^p : \Lambda ^\infty \to
\Lambda ^\infty $, by
  $$
  \sigma ^p (x)(m, n) = x(p + m, p + n), \for x\in \Lambda ^\infty , \for (m, n)\in \Omega _k.
  $$

In what follows we would like to relate the path space $\Lambda ^\infty $ to an ${\bf N}^k$-subshift.  In order to do so
we begin by introducing the notation
  $
  \1:=(1,1,\ldots ,1),
  $
  which we use in defining our alphabet $\Sigma $, by
  $$
  \Sigma =\big \{\lambda \in \Lambda : d(\lambda ) = \1\big \}.
  $$

Given any path $x\in \Lambda ^\infty $, we may consider the element $\xi _x$ of $\Sigma ^{{\bf N}^k}$, defined by
  $$
  \xi _x(n)=x(n, n{+}\1), \for n\in {\bf N}^k,
  $$

\beginmypix x from 0 to 4 y from -1 to 4.6 xscale 1.5truecm yscale 1.5truecm \put{$\bullet$} at 0 0 \plot 0 0 1 0 /
\plot 0 0 0 1 / \put{${\scriptstyle \xi _x(0, 0)}$} at 0.5 0.5 \put{$\bullet$} at 0 1 \plot 0 1 1 1 / \plot 0 1 0 2 /
\put{${\scriptstyle \xi _x(0, 1)}$} at 0.5 1.5 \put{$\bullet$} at 0 2 \plot 0 2 1 2 / \plot 0 2 0 3 /
\put{${\scriptstyle \xi _x(0, 2)}$} at 0.5 2.5 \put{$\bullet$} at 0 3 \plot 0 3 1 3 / \plot 0 3 0 4 /
\put{${\scriptstyle \xi _x(0, 3)}$} at 0.5 3.5 \put{$\bullet$} at 0 4 \plot 0 4 1 4 / \plot 0 4 0 4.3 / \put{$\bullet$}
at 1 0 \plot 1 0 2 0 / \plot 1 0 1 1 / \put{${\scriptstyle \xi _x(1, 0)}$} at 1.5 0.5 \put{$\bullet$} at 1 1 \plot 1 1 2
1 / \plot 1 1 1 2 / \put{${\scriptstyle \xi _x(1, 1)}$} at 1.5 1.5 \put{$\bullet$} at 1 2 \plot 1 2 2 2 / \plot 1 2 1 3
/ \put{${\scriptstyle \xi _x(1, 2)}$} at 1.5 2.5 \put{$\bullet$} at 1 3 \plot 1 3 2 3 / \plot 1 3 1 4 /
\put{${\scriptstyle \xi _x(1, 3)}$} at 1.5 3.5 \put{$\bullet$} at 1 4 \plot 1 4 2 4 / \plot 1 4 1 4.3 / \put{$\bullet$}
at 2 0 \plot 2 0 3 0 / \plot 2 0 2 1 / \put{${\scriptstyle \xi _x(2, 0)}$} at 2.5 0.5 \put{$\bullet$} at 2 1 \plot 2 1 3
1 / \plot 2 1 2 2 / \put{${\scriptstyle \xi _x(2, 1)}$} at 2.5 1.5 \put{$\bullet$} at 2 2 \plot 2 2 3 2 / \plot 2 2 2 3
/ \put{${\scriptstyle \xi _x(2, 2)}$} at 2.5 2.5 \put{$\bullet$} at 2 3 \plot 2 3 3 3 / \plot 2 3 2 4 /
\put{${\scriptstyle \xi _x(2, 3)}$} at 2.5 3.5 \put{$\bullet$} at 2 4 \plot 2 4 3 4 / \plot 2 4 2 4.3 / \put{$\bullet$}
at 3 0 \plot 3 0 4 0 / \plot 3 0 3 1 / \put{${\scriptstyle \xi _x(3, 0)}$} at 3.5 0.5 \put{$\bullet$} at 3 1 \plot 3 1 4
1 / \plot 3 1 3 2 / \put{${\scriptstyle \xi _x(3, 1)}$} at 3.5 1.5 \put{$\bullet$} at 3 2 \plot 3 2 4 2 / \plot 3 2 3 3
/ \put{${\scriptstyle \xi _x(3, 2)}$} at 3.5 2.5 \put{$\bullet$} at 3 3 \plot 3 3 4 3 / \plot 3 3 3 4 /
\put{${\scriptstyle \xi _x(3, 3)}$} at 3.5 3.5 \put{$\bullet$} at 3 4 \plot 3 4 4 4 / \plot 3 4 3 4.3 / \put{$\bullet$}
at 4 0 \plot 4 0 4.3 0 / \plot 4 0 4 1 / \put{$\bullet$} at 4 1 \plot 4 1 4.3 1 / \plot 4 1 4 2 / \put{$\bullet$} at 4 2
\plot 4 2 4.3 2 / \plot 4 2 4 3 / \put{$\bullet$} at 4 3 \plot 4 3 4.3 3 / \plot 4 3 4 4 / \put{$\bullet$} at 4 4 \plot
4 4 4.3 4 / \plot 4 4 4 4.3 / \put{\eightrm A representation of the path $\xi _x$} at 2 -0.4 \endmypix

For every $p$ in ${\bf N}^k$, and every $x$ in $\Lambda ^\infty $, one easily checks that
  $$
  \theta _p(\xi _x) = \xi _{\sigma ^p(x)},
  \equationmark Covar
  $$
  from where one deduces that the subset of $\Sigma ^{{\bf N}^k}$, given by
  $$
  X_\Lambda =\big \{\xi _x:x\in \Lambda ^\infty \big \}
  $$
  is invariant under the Bernoulli action of ${\bf N}^k$.

\state Proposition \label XLambda
  The set $X_\Lambda $ introduced above is an ${\bf N}^k$-subshift and the correspondence
  $$
  \Xi : x\in \Lambda ^\infty \mapsto \xi _x\in X_\Lambda
  $$
  is a homeomorphism.  In addition, $\Xi $ is covariant relative to the action $\sigma $ on $\Lambda ^\infty $ and the
Bernoulli action on $X_\Lambda $.

\Proof
  Having already observed that $X_\Lambda $ is invariant under the Bernoulli action, we will next prove that $X_\Lambda
$ is closed.  We then pick any $\xi $ in the closure of $X_\Lambda $, so there exists a sequence $\{x_i\}_i$ in $\Lambda
^\infty $ such that $\{\xi _{x_i}\}_i$ converges to $\xi $.

  For each $j$ in ${\bf N}$, write $j\1$ for the element $(j,j,\ldots ,j)$ of ${\bf N}^k$, and put
  $
  \lambda _j=\xi (j\1).
  $
  We then claim that $s(\lambda _j)=r(\lambda _{j+1})$, for every $j$, where $s$ and $r$ refer to the source and range
maps relative to the category $\Lambda $.  To see this, let $i_0$ be large enough, so that
  $$
  \xi _{x_i}(j\1)=\xi (j\1), \and \xi _{x_i}(j\1{+}\1)=\xi (j\1{+}\1),
  $$
  for every $i\geq i_0$.  Choosing any $i\geq i_0$, we then have that
  $$
  s(\lambda _j)=
  s\big(\xi (j\1)\big) =
  s\big(\xi _{x_i}(j\1)\big) =
  s\big({x_i}(j\1, j\1{+}\1)\big) = \cdots
  $$
  Recall that $x_i$ is a path, hence a functor from $\Omega _k$ to $\Lambda $.  Since the morphisms $(j\1, j\1{+}\1)$
and $(j\1{+}\1, j\1+2\1)$ may be composed in $\Omega _k$, we have that $x_i(j\1, j\1{+}\1)$ and $x_i(j\1{+}\1, j\1+2\1)$
may be composed in $\Lambda $, so the source of the former must coincide with the range of the latter, whence the above
equals
  $$
  \cdots = r\big({x_i}(j\1{+}\1, j\1+2\1)\big) =
  r\big(\xi _{x_i}(j\1{+}\1)\big) =
  r\big(\xi (j\1{+}\1)\big) =
  r(\lambda _j),
  $$
  thus proving our claim.  By the paragraph after \cite [Remarks 2.2]{KP}, we conclude that there exists a path $x$ in
$\Lambda ^\infty $ such that
  $x(j\1, j\1{+}\1)=\lambda _j$, for all $j$, and we next claim that $\xi _x=\xi $.  Given any $n$ in ${\bf N}^k$,
choose some $j$ in ${\bf N}$ such that $n\leq j\1$, and observe that
  $$
  \eqnarray={
  \lambda _0\lambda _1\ldots \lambda _j
  & x(0\1,\1)\ x(\1, 2\1)\ \cdots \ x(j\1,j\1{+}\1) \cr
  & x(0\1,j\1{+}\1) \cr
  & x(0\1,n)\ x(n, n{+}\1)\ x(n{+}\1, j\1{+}\1) \cr
  & x(0\1,n)\ \xi _x(n)\ x(n{+}\1, j\1{+}\1). \cr}
  $$
  We next choose $i_0$ large enough, so that $\xi _{x_i}$ coincides with $\xi $ on
  $\{0\1, \1, \ldots , j\1,n \}$, for every $i\geq i_0$.  Therefore
  $$
  \eqnarray={
  \lambda _0\lambda _1\ldots \lambda _j
  & \xi (0\1)\xi (\1)\ldots \xi (j\1) \cr
  & \xi _{x_i}(0\1)\ \xi _{x_i}(\1)\ldots \xi _{x_i}(j\1) \cr
  & x_i(0\1, \1)\ x_i(\1, 2\1) \ldots x_i(j\1, j\1{+}\1) \cr
  & x_i(0\1,j\1{+}\1) \cr
  & x_i(0\1,n)\ x_i(n, n{+}\1)\ x_i(n{+}\1, j\1{+}\1). \cr}
  $$
  Contrasting our last two calculations, and invoking the unique factorization property, we deduce that
  $$
  \xi _x(n)=x_i(n, n{+}\1) = \xi _{x_i}(n)=\xi (n).
  $$
  This concludes the proof of the claim according to which $\xi _x=\xi $, whence $\xi $ lies in $X_\Lambda $, and so we
see that $X_\Lambda $ is closed.

The last sentence of the statement has already been verified in \ref {Covar}, so we are finally left with the task of
proving our correspondence $\Xi $ to be a homeomorphism.  In order to do this we first observe that $\Xi $ is injective
by the paragraph after \cite [Remarks 2.2]{KP}.  We will next prove that $\Xi $ is an open mapping.  For this, recall
from \cite [Definitions 2.4]{KP} that the topology of $\Lambda ^\infty $ is generated by the \emph {cylinders},
  $$
  Z(\lambda ) = \{x\in \Lambda ^\infty : x\big(0, d(\lambda )\big) = \lambda \},
  $$
  as $\lambda $ range in $\Lambda $.  Given any $y$ in $\Lambda ^\infty $, and given any open set $U\subseteq \Lambda
^\infty $ containing $y$, one may then find some $\lambda $ such that
  $$
  y\in Z(\lambda )\subseteq U.
  $$
  Choosing $j$ such that $j\1\geq d(\lambda )$, and setting $\mu =y(0,j\1)$, we then have that
  $$
  \mu = y(0,j\1) = y\big(0,d(\lambda )\big)\ y\big(d(\lambda ), j\1\big) = \lambda \, y\big(d(\lambda ), j\1\big),
  $$
  from where we see that
  $$
  y\in Z(\mu )\subseteq Z(\lambda ),
  $$
  and we conclude that the cylinders of the form $Z(\mu )$, with $d(\mu )$ a multiple of $\1$, also form a basis for the
topology of $\Lambda ^\infty $.  In order to prove that $\Xi $ is an open map, it is therefore enough to verify that
$\Xi \big(Z(\mu )\big)$ is open for every such $\mu $.

Given $\mu $ in $\Lambda $ as above, i.e.~with $d(\mu )=j\1$, use the factorization property to write
  $$
  \mu =\lambda _0\lambda _1\ldots \lambda _{j-1},
  $$
  with $d(\lambda _i)=\1$, and observe that
  $$
  \eqnarray\Leftrightarrow{
  x\in Z(\mu )
  & x(0, j\1)= \lambda _0\lambda _1\ldots \lambda _{j-1} \cr
  & x(i\1,i\1{+}\1)=\lambda _i, \for i=0,\cdots , j{-}1\phantom {.} \cr
  & \xi _x(i\1)=\lambda _i, \hfill \for i=0,\cdots , j{-}1. \cr
  }
  $$
  Setting
  $$
  V=\big \{\xi \in \Sigma ^{{\bf N}^k}: \xi (i\1)=\lambda _i,\ \forall i=0,\cdots , j{-}1\big \},
  $$
  which is clearly open in $\Sigma ^{{\bf N}^k}$, we then deduce that $x\in Z(\mu )$ if and only if $\xi _x\in V$. It
then follows that $\Xi \big(Z(\mu )\big) = X_\Lambda \cap V$, proving that $\Xi $ is an open mapping, as claimed.

Finally, leaving for the reader the easy task of verifying that $\Xi $ is continuous, the proof is concluded.  \endProof

Having twice resorted to the paragraph after \cite [Remarks 2.2]{KP}, we nave not yet exhausted its consequences from
our point of view.  A further, and major consequence is the content of our next result.

Identifying the monoid ${\bf N}$ as a submonoid of ${\bf N}^k$ via the correspondence $i \leftrightarrow i\1$, we will
shortly refer to the restriction map $\rho \sub {\bf N}$, introduced in \ref {RestrMap}.

\state Proposition \label HighCore
  Let $A=\{A_{\lambda \mu }\}_{\lambda ,\mu \in \Sigma }$ be the matrix given by
  $$
  A_{\lambda \mu } =
  \clauses {
    \cl 1 if {s(\lambda )=r(\mu ),}
    \cl 0 otherwise, {}
    }
  $$
  and let $X_A\subseteq \Sigma ^{\bf N}$ be the Markov space for $A$.  Then $\rho \sub {\bf N}(X_\Lambda )=X_A$, and
$\rho \sub {\bf N}$ is a homeomorphism from $X_\Lambda $ onto $X_A$.

\Proof
  We first observe that, for every $x$ in $\Lambda ^\infty $, one has
  $$
  \rho \sub {\bf N}(\xi _x) = \big (x(i\1,i\1{+}\1)\big )_{i\in {\bf N}},
  $$
  which is evidently in $X_A$, so we see that $\rho \sub {\bf N}(X_\Lambda )\subseteq X_A$.

Given any $\lambda =(\lambda _i)_{i\in {\bf N}}\in X_A$, the paragraph after \cite [Remarks 2.2]{KP} produces a path $x$
in $\Lambda ^\infty $ such that $x(i\1,i\1{+}\1)=\lambda _i$, for all $i$ in ${\bf N}$, so that $\rho \sub {\bf N}(\xi
_x) = \lambda $, which in turn proves that $\rho \sub {\bf N}(X_\Lambda )=X_A$.

The uniqueness of the path $x$ obtained above implies that $\rho \sub {\bf N}$ is one-to-one, so it remains to prove
that $\rho \sub {\bf N}$ is a homeomorphism.  Continuity not being an issue, we focus on proving continuity of the
inverse map.  In order to do so let us review the above construction of $x $ from any given $\lambda $ in $X_A$: given
any $(m, n)\in \Omega _k$, one chooses $j$ such that $j\1\geq n$, and uses the unique factorization property to write
  $$
  \lambda _0\lambda _1\ldots \lambda _{j-1} = \mu x(m, n)\nu ,
  $$
  with $d(\mu )=m$, \ $d\big(x(m, n)\big) = n-m$, and $d(\nu ) = j\1-n$.  The resulting map $x:\Omega _k\to \Lambda $ is
then a path satisfying $\rho \sub {\bf N}(\xi _x) = \lambda $, whence $\xi _x=\rho \sub {\bf N}\inv (\lambda )$, and we
must then prove that $\xi _x$ varies continuously with $\lambda $.  Since $\xi _x$ lives in the product space $\Sigma
^{{\bf N}^k}$, all we need to do is show that $\xi _x(n)$ is continuous as a function of $\lambda $, for every $n$ in
${\bf N}^k$.  Recalling that $\xi _x(n)=x(n, n{+}\1)$, notice that the above recipe to construct $x(n, n{+}\1)$ depends
only on $(\lambda _0,\lambda _1,\ldots , \lambda _{j-1})$, where $j$ is any integer such that $j\1\geq n{+}\1$.  The
function producing $\xi _x(n)$ from $\lambda $ therefore factors as the composition
  $$
  X_A\to \Sigma ^j \to \Sigma ,
  $$
  where the leftmost arrow is the projection on the first $j$ coordinates while the rightmost one corresponds to the
above recipe producing $x(n, n{+}\1)$ from $(\lambda _0,\lambda _1,\ldots , \lambda _j)$.  Since the projection is
continuous and $\Sigma ^j$ carries the discrete topology, continuity follows.  \endProof

\section Higher rank subshifts

Motivated by the example of the ${\bf N}^k$-subshift discussed in the previous section, we make the following:

\definition \label DefRankShift
  Let $k$ be a positive integer.  By a \emph {subshift of rank $k$}, or a \emph {$k$-subshift}, we shall mean a
  ${\bf N}^k$-subshift $X$ on a given alphabet $\Sigma $ such that
  \izitem
  \zitem the restriction map $\rho \sub {\bf N}$ is injective on $X$,
  \zitem $\rho \sub {\bf N}(X)$ is closed in $\Sigma ^{\bf N}$, and
  \zitem $\rho \sub {\bf N}$ is a homeomorphism from $X$ to $\rho \sub {\bf N}(X)$.

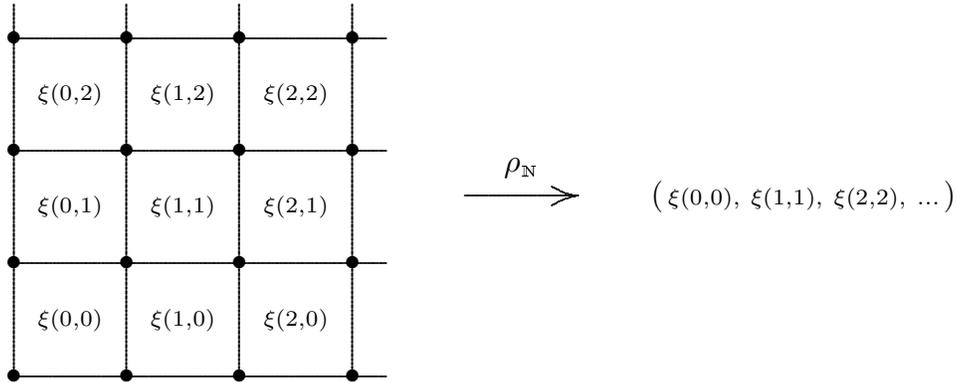
\beginmypix x from 0 to 8 y from -0.5 to 3.5 xscale 1.5truecm yscale 1.5truecm \put{$\bullet$} at 0 0 \plot 0 0 1 0 /
\plot 0 0 0 1 / \put{${\scriptstyle \xi (0, 0)}$} at 0.5 0.5 \put{$\bullet$} at 0 1 \plot 0 1 1 1 / \plot 0 1 0 2 /
\put{${\scriptstyle \xi (0, 1)}$} at 0.5 1.5 \put{$\bullet$} at 0 2 \plot 0 2 1 2 / \plot 0 2 0 3 / \put{${\scriptstyle
\xi (0, 2)}$} at 0.5 2.5 \put{$\bullet$} at 0 3 \plot 0 3 1 3 / \plot 0 3 0 3.3 / \put{$\bullet$} at 1 0 \plot 1 0 2 0 /
\plot 1 0 1 1 / \put{${\scriptstyle \xi (1, 0)}$} at 1.5 0.5 \put{$\bullet$} at 1 1 \plot 1 1 2 1 / \plot 1 1 1 2 /
\put{${\scriptstyle \xi (1, 1)}$} at 1.5 1.5 \put{$\bullet$} at 1 2 \plot 1 2 2 2 / \plot 1 2 1 3 / \put{${\scriptstyle
\xi (1, 2)}$} at 1.5 2.5 \put{$\bullet$} at 1 3 \plot 1 3 2 3 / \plot 1 3 1 3.3 / \put{$\bullet$} at 2 0 \plot 2 0 3 0 /
\plot 2 0 2 1 / \put{${\scriptstyle \xi (2, 0)}$} at 2.5 0.5 \put{$\bullet$} at 2 1 \plot 2 1 3 1 / \plot 2 1 2 2 /
\put{${\scriptstyle \xi (2, 1)}$} at 2.5 1.5 \put{$\bullet$} at 2 2 \plot 2 2 3 2 / \plot 2 2 2 3 / \put{${\scriptstyle
\xi (2, 2)}$} at 2.5 2.5 \put{$\bullet$} at 2 3 \plot 2 3 3 3 / \plot 2 3 2 3.3 / \put{$\bullet$} at 3 0 \plot 3 0 3.3 0
/ \plot 3 0 3 1 / \put{$\bullet$} at 3 1 \plot 3 1 3.3 1 / \plot 3 1 3 2 / \put{$\bullet$} at 3 2 \plot 3 2 3.3 2 /
\plot 3 2 3 3 / \put{$\bullet$} at 3 3 \plot 3 3 3.3 3 / \plot 3 3 3 3.3 / \arrow <0.3cm> [0.25, 0.75] from 4 1.6 to 5
1.6 \put{$\rho \sub {\bf N}$} at 4.5 1.85 \put{($\scriptstyle \, \xi (0, 0), \ \xi (1, 1), \ \xi (2, 2), \ \ldots\, $)}
at 7 1.6 \endmypix

Since Markov spaces are automatically closed, we conclude from \ref {HighCore} that the ${\bf N}^k$-subshift $X_\Lambda
$ built from a $k$-graph $\Lambda $ is an example of a $k$-subshift.

Observe also that, in case the alphabet $\Sigma $ is finite, then any ${\bf N}^k$-subshift is compact, hence \ref
{DefRankShift.ii} is automatically true, while \ref {DefRankShift.iii} follows from \ref {DefRankShift.i}.  In other
words, when the alphabet is finite, \ref {DefRankShift.ii-iii} could be omitted from the above definition without any
consequences.

\fix Let us now fix a subshift $X$ of rank $k$ on the alphabet $\Sigma $.

\medskip \noindent Setting $Y=\rho \sub {\bf N}(X)$, we have by \ref {RestrictPShift.i} that $Y$ is invariant under the
Bernoulli action of $\N $, that is, invariant under the usual shift map
  $$
  S:\Sigma ^{\bf N}\to \Sigma ^{\bf N},
  $$
  and since $Y$ is also closed by assumption, we have that $Y$ is a ${\bf N}$-subshift, that is, a classical subshift.

{\it A priori}, it does not make sense to ask whether or not $\rho \sub {\bf N}$ is covariant, since the monoids acting
on the Bernoulli spaces $\Sigma ^{{\bf N}^k}$ and $\Sigma ^{\bf N}$ are not the same.  But if we consider only the
smaller monoid, namely ${\bf N}$, then covariance clearly holds, and in particular
  $$
  \rho \sub {\bf N}\theta _\1 = S \rho \sub {\bf N},
  \equationmark CovarRestr
  $$
  as the reader may easily verify.

Furthermore, since $\rho \sub {\bf N}$ is a homeomorphism from $X$ to $Y$, the Bernoulli action of ${\bf N}^k$ on $X$
may be transfered to $Y$, so we get an action $\tau $ of ${\bf N}^k$ on $Y$ such that the diagram

\beginmypix x from -200 to 400 y from -100 to 350 xscale 0.012truecm yscale -0.012truecm \map $X$ 000 000 $Y$ 200 000
\lbl <0pt,10pt> {$\rho \sub {\bf N}$} \map $X$ 000 200 $Y$ 200 200 \lbl <0pt,-10pt> {$\rho \sub {\bf N}$} \flexa 000 000
000 200 \lbl <-10pt,0pt> {$\theta _n$} \flexa 200 000 200 200 \lbl <10pt,0pt> {$\tau _ n$} \advseqnumbering \put
{\eightrm Diagram (\current)} at 100 290 \label Diagrama \endmypix

\noindent commutes for every $n$ in ${\bf N}^k$.

\state Theorem \label Factorize
  Let $X$ be a $k$-subshift on the alphabet $\Sigma $, and put $Y=\rho \sub {\bf N}(X)$.  Then:
  \izitem
  \zitem $Y$ is a classical subshift.
  \zitem For any integer $i$ with $1\leq i\leq k$, let $e_i$ be $i^{th}$ canonical basis vector of ${\bf N}^k$, and put
$S_i=\tau _{e_i}$.  Then the $S_i$ are pairwise commuting, continuous maps from $Y$ to $Y$, and
  $$
  S_1S_2\cdots S_k=S,
  $$
  where $S$ is the restriction of the shift map to $Y$ (here denoted simply by $S$, by abuse of language).

\Proof The first point is an obvious consequence of the definitions and of \ref {RestrictPShift.i}.  It was included
here only for future reference.  Regarding (ii) we have
  $$
  S_1S_2\cdots S_k=\tau (e_1)\cdots \tau (e_k)=\tau (e_1+\cdots +e_k)=\tau (\1)=\rho \sub {\bf N}\theta _\1\rho \sub
{\bf N}\inv \={CovarRestr} S.
  \closeProof
  $$
  \endProof

Given that the $S_i$ commute among themselves, it follows that $S_i$ also commutes with $S$, so this is turns out to be
strongly related to \ref{Commutativity}, and hence also to cellular automata by \ref{Hedlund}.

\section Cellular automaton factorization of Markov subshifs associated to $k$-graphs

\label Factoriza We have already mentioned that the space $X_\Lambda $ built from a $k$-graph $\Lambda $ is a subshift
of rank $k$.  Moreover, by \ref {HighCore}, we have that $\rho \sub {\bf N}(X_\Lambda )$ is a Markov space.  We may then
use \ref {Factorize} to produce a factorization of Markov's shift and we will now show that each $S_i$ occuring in \ref
{Factorize} is in fact a cellular automaton with anticipation $1$.
  This could be obtained by \ref{Hedlund}, but we may in fact produce the block codes directly.

We first observe that the language of a Markov subshift, such as $\rho \sub {\bf N}(X_\Lambda )$, is governed by its
matrix, and in particular
  $$
  \L_1\big(\rho (X_\Lambda )\big) = \{(\lambda _0,\lambda _1) \in \Sigma ^2: s(\lambda _0)=r(\lambda _1)\}.
  $$

We therefore define, for each $i=1,\ldots , k$,
  $$
  \varphi _i:\L_1\big(\rho (X_\Lambda )\big)\to \Sigma
  $$
  as follows: given $(\lambda ,\mu )$ in $\L_1\big(\rho (X_\Lambda )\big)$, we have that $\lambda \mu \in \Lambda $, and
clearly $d(\lambda \mu )=\1+\1$.  Writing
  $$
  \1+\1 = e_i + \1 + (\1 -e_i),
  $$
  the unique factorization property allows us to write
  $
  \lambda \mu = \alpha \beta \gamma ,
  $
  where $\alpha $, $\beta $ and $\gamma $ lie in $\Lambda $, $d(\alpha )=e_i$, $d(\beta )= \1$, and $d(\gamma )=
\1-e_i$.  We then set
  $$
  \varphi _i(\lambda ,\mu )=\beta .
  $$

We will now show that each $S_i$ coincides with the cellular automaton relative to the sliding block code $\varphi _i$.
In order to do this, observe that each $S_i$ is officially defined as
  $$
  S_i = \tau _{e_i} = \rho \sub {\bf N}\theta _{e_i}\rho \sub {\bf N}\inv .
  $$
  Given any $y$ in $\rho \sub {\bf N}(X_\Lambda )$, we may write $y=\rho \sub {\bf N}(\xi )$, for some $\xi $ in
$X_\Lambda $, and we may further write $\xi =\xi _x$, for some $x\in \Lambda ^\infty $.  In other words, $y=\rho \sub
{\bf N}(\xi _x)$.  We then have for every $j$ in ${\bf N}$, that
  $$
  S_i(y)\calcat j =
  \rho \sub {\bf N}\theta _{e_i}\rho \sub {\bf N}\inv (y)\calcat j =
  \theta _{e_i}\rho \sub {\bf N}\inv (y)\calcat {j\1} =
  \rho \sub {\bf N}\inv (y)\calcat {j\1+e_i} \quebra=
  \xi _x(j\1+e_i) =
  x(j\1+e_i, j\1+e_i{+}\1).
  \equationmark SiY
  $$
  On the other hand, notice that
  $$
  y(j)\, y(j{+}1) =
  \xi _x(j\1)\, \xi _x(j\1{+}\1) =
  x(j\1,j\1{+}\1)\,x(j\1{+}\1,j\1{+}2\1) =
  x(j\1, j\1{+}2\1) \quebra =
  x(j\1, j\1{+}e_i)\,x\big(j\1{+}e_i, j\1{+}e_i{+}\1)\,x(j\1{+}e_i{+}\1,j\1{+}2\1).
  $$
  By the unique factorization property we then have that
  $$
  \varphi _i\big(y(j),y(j{+}1)\big) = x(j\1{+}e_i, j\1{+}e_i{+}\1) \={SiY} S_i(y)\calcat j,
  $$
  thus proving that indeed $S_i$ is the cellular automaton associated to the block code $\varphi _i$, as claimed.

\section Constructing $k$-subshifts from factorizations of the shift

Motivated by \ref {Factorize} we introduce the following concept:

\definition \label Kautomaton
  A \emph {\wkautomat on} over a given alphabet $\Sigma $ is a $k{+}1$-tuple
  $$
  (Y;\ S_1,S_2,\ldots ,S_k),
  $$
  where $Y$ is a closed subset of $\Sigma ^{\bf N}$, invariant under the shift $S$ (i.e., $Y$ is a classical subshift),
and the $S_i$ are pairwise commuting continuous maps from $Y$ to $Y$ such that
  $$
  S_1S_2\cdots S_k=S.
  $$

As seen in \ref {Factorize}, every $k$-subshift leads to a \wkautomat on, and it is our plan to show that $k$-subshifts
are essentially the same thing as \wkautomat a. As a first step let us show how to construct a $k$-subshift given a
\wkautomat a $(Y;\ S_1,S_2,\ldots ,S_k)$, which we consider fixed for the time being.

For each $n=(n_1,n_2,\ldots ,n_k)\in {\bf N}^k$, put
  $$
  \alpha _n=S_1^{n_1}S_2^{n_2}\cdots S_k^{n_k},
  $$
  so that $\alpha $ is an action of ${\bf N}^k$ on $Y$.  We next let
  $
  \Psi :Y\to \Sigma ^{{\bf N}^k}
  $
  be defined by
  $$
  \Psi (y)\calcat n = \alpha _n(y)\calcat 0, \for y \in Y, \for n \in {\bf N}^k.
  $$

\state Proposition \label ConstructKSubshift
  Setting $X=\Psi (Y)$, one has that $X$ is a subshift of rank $k$.  In addition we have that $\rho \sub {\bf N}(X)=Y$
and, regarding the canonical action $\tau $ introduced in Diagram \ref {Diagrama}, one has that $\tau _{e_i}=S_i$, for
every $i=1, \ldots , k$.

\Proof
  Given $y$ in $Y$, and given $p$ and $t$ in ${\bf N}^k$, observe that
  $$
  \theta _p\big(\Psi (y)\big)\calcat t = \Psi (y)\calcat {t+p} = \alpha _{t+p}(y)\calcat 0 = \alpha _t\big(\alpha
_p(y)\big)\calcat 0 =
  \Psi \big(\alpha _p(y)\big)\calcat t,
  \equationmark CovarTau
  $$
  so we deduce that $\theta _p\big(\Psi (y)\big)= \Psi \big(\alpha _p(y)\big)$, from where we see that the range of
$\Psi $, also known as $X$, is invariant under the Bernoulli action of ${\bf N}^k$.

In order to prove that $\rho \sub {\bf N}(X)=Y$, it suffices to verify that
  $$
  \rho \sub {\bf N}(\Psi (y))=y, \for y\in Y,
  \equationmark SectRho
  $$
  but this follows from the following computation, where $j\in {\bf N}$:
  $$
  \rho \sub {\bf N}(\Psi (y))\calcat j = \Psi (y)\calcat {j\1} =
  \alpha _{j\1}(y)\calcat 0 =
  S_1^jS_2^j\cdots S_k^j(y)\calcat 0 =
  S^j(y)\calcat 0 =
  y(j).
  $$

In order to prove that $X$ is closed, suppose that $\{y_i\}_i$ is a sequence in $Y$ such that $\{\Psi (y_i)\}_i$
converges to some $x$ in $\Sigma ^{{\bf N}^k}$.  Then
  $$
  y:=\rho \sub {\bf N}(x) = \lim _i \rho \sub {\bf N}(\Psi (y_i)) \={SectRho}
  \lim _i y_i,
  $$
  so $y\in Y$, and we claim that $x=\Psi (y)$.  In fact, for every $n$ in ${\bf N}^k$, we have that
  $$
  x(n) =
  \lim _i \Psi (y_i)\calcat n =
  \lim _i \alpha _n(y_i)\calcat 0 =
  \alpha _n(y)\calcat 0 =
  \Psi (y)\calcat n,
  $$
  proving the claim.

So far we have thus proven that $X$ is a ${\bf N}^k$-subshift.  In order to show that it is a subshift of rank $k$, we
must verify \ref {DefRankShift.i--iii}.  By \ref {SectRho} we see that $\Psi $ is one-to-one, and since it is onto $X$,
by definition, it follows that $\Psi $ is bijective.  Employing \ref {SectRho} once more we deduce that $\rho \sub {\bf
N}|_X=\Psi \inv $, from where \ref {DefRankShift.i} follows.

Noticing that both $\rho \sub {\bf N}$ and $\Psi $ are clearly continuous, we obtain \ref {DefRankShift.iii}, while \ref
{DefRankShift.ii} follows from the facts that $\rho \sub {\bf N}(X)=Y$, and that $Y$ is closed by assumption.

The last part of the statement may now be proved as follows: for $y\in Y$, and $j\in {\bf N}$, one has
  $$
  \tau _{e_i}(y)\calcat j =
  \rho \sub {\bf N}\theta _{e_i}\rho \sub {\bf N}\inv (y)\calcat j =
  \theta _{e_i}\rho \sub {\bf N}\inv (y)\calcat {j\1} =
  \theta _{e_i}\Psi (y)\calcat {j\1} \={CovarTau}
  \Psi \big(\alpha _{e_i}(y)\big)\calcat {j\1} \quebra =
  \Psi \big(S_i(y)\big)\calcat {j\1} =
  \alpha _{j\1}(S_i(y))\calcat 0 =
  S^j\big(S_i(y)\big)\calcat 0 =
  S_i(y)\calcat j,
  $$
  so $\tau _{e_i}= S_i$.  \endProof

\state Proposition
  Suppose that $X$ and $X'$ are $k$-subshifts over the same alphabet $\Sigma $ such that $\rho \sub {\bf N}(X)=\rho \sub
{\bf N}(X')$, and such that the canonical actions $\tau $ and $\tau '$ coincide.  Then $X=X'$.

\Proof
  Given $x\in X$, and $n\in {\bf N}^k$, we have seen in \ref {Diagrama} that
  $
  \tau _n \rho \sub {\bf N}(x) = \rho \sub {\bf N}\theta _n (x).
  $
  Therefore
  $$
  x(n) = \theta _n (x)\calcat {0\1} = \rho \sub {\bf N}\theta _n (x)\calcat 0 =
  \tau _n \rho \sub {\bf N}(x)\calcat 0.
  $$
  This says that $x$ may be reconstructed from $\rho \sub {\bf N}(x)$ together with the canonical action $\tau $, from
where the result follows.  \endProof

Summarizing our last two results we have:

\state Corollary \label SameThing
  Given $k\geq 1$ and an alphabet $\Sigma $, the corresponence
  $$
  X \mapsto \big(\rho \sub {\bf N}(X);\ \tau _{e_1},\tau _{e_2},\ldots ,\tau _{e_k}\big)
  $$
  establishes a one-to-one correspondence from the collection of all $k$-subshifts $X\subseteq \Sigma ^{{\bf N}^k}$ onto
the collection of all \wkautomat a on the alphabet $\Sigma $.

The adjective ``weak'' employed in Definition \ref{Kautomaton} is meant to highlight the fact that the $S_i$ mentioned
there are not actually cellular automata, although they share with the latter the important property of commuting with
the shift, a property we saw in \ref{Hedlund} to characterize true cellular automata in the uniformly continuous case.

Nevertheless it is interesting to determine necessary and sufficient conditions on a given $k$-subshift for the maps
$S_i$ in the \wkautomat a associated to it by \ref{SameThing} to be actual cellular automata.
  In order to do this we must first consider a metric on $\Sigma ^{{\bf N}^k}$ as follows.

\definition
  We shall say that two elements $x$ and $y$ of\/ $\Sigma ^{{\bf N}^k}$ \emph{agree} on a given subset $A\subseteq {\bf
N}^k$ , whenever
  $x(n)=y(n)$, for all $n$ in $A$.
  If $p$ is the largest integer such that $x$ and $y$ agree on the subset
  $$
  B_p := \big\{n=(n_1, n_2, \ldots , n_k)\in {\bf N}^k: n_i\leq p, \hbox{ for all }i\big\},
  $$
  we put
  $$
  d(x,y) = 2^{-p},
  $$
  with the understanding that if $x=y$, then $p=\infty $, in which case $d(x,y) = 0$.

As in the case of section \ref{CellAutomataSection}, one proves that $d$ is a metric on
  $\Sigma ^{{\bf N}^k}$, compatible with the product topology, and we may then speak of uniformly continuous functions
defined, or taking values, in $\Sigma ^{{\bf N}^k}$.

The above choice of the $B_p$ is not so crucial except for the fact that the $B_p$ form an increasing sequence of finite
subsets of ${\bf N}^k$, whose union coincides with ${\bf N}^k$.  Any other choice of finite subsets with these
properties may also be used to define a metric on $\Sigma ^{{\bf N}^k}$, which in turn induce the same uniform structure
\cite[Chapter 6]{Kelley} on $\Sigma ^{{\bf N}^k}$.  The common underlying uniform structure is in fact what really
matters here, and it says that two points $x$ and $y$ are close if and only if they agree on a large finite subset of
${\bf N}^k$.

We shall however not make any explicit use of uniform structures in this work, beyond the elementary observation that
  $$
  d(x,y)\leq 2^{-p} \iff \hbox{ $x$ and $y$ agree on $B_p$},
  \equationmark BabyUniformStru
  $$
  for every $x$ and $y$ in $\Sigma ^{{\bf N}^k}$.

\state Proposition
  The restriction map
  $$
  \rho \sub{\bf N}: \Sigma ^{{\bf N}^k} \to \Sigma ^{\bf N},
  $$
  defined in \ref{RestrMap} is uniformly continuous.

\Proof Given $x$ and $y$ in $\Sigma ^{{\bf N}^k}$, let $p$ be the largest integer such that $x$ and $y$ agree on $B_p$.
Therefore $\rho \sub{\bf N}(x)$ and $\rho \sub{\bf N}(y)$ obviously agree on $\{0, 1, \ldots , p\}$, so
  $$
  d\big(\rho \sub{\bf N}(x),\rho \sub{\bf N}(y)\big)\leq 2^{-p} = d(x, y).
  $$
  This proves that $\rho \sub{\bf N}$ is in fact contractive, hence uniformly continuous.
   \endProof

If $X\subseteq \Sigma ^{{\bf N}^k}$ is a $k$-subshift, then $\rho \sub{\bf N}$ is a homeomorphism on $X$ by definition,
so $\rho \sub{\bf N}\inv$ is continuous on $\rho \sub{\bf N}(X)$, although perhaps not uniformly.

\definition \label DefineGoodStuff
  \iaitem
  \aitem By a \emph{uniform $k$-subshift} we shall mean a $k$-subshift $X\subseteq \Sigma ^{{\bf N}^k}$ such that the
inverse of $\rho \sub{\bf N}$ is uniformly continuous on $\rho \sub{\bf N}(X)$.
  \aitem By a \emph{$k$-automaton} we shall mean a \wkautomat on $(Y;S_1,S_2,\ldots ,S_k)$ such that each $S_i$ is
actually a cellular automaton.

In what follows we will show that the two concepts just defined are related to each other by the same process involved
in \ref {SameThing}.

\state Proposition \label CharacStrongAutomata
  Let $X\subseteq \Sigma ^{{\bf N}^k}$ be a $k$-subshift.  Then $X$ is a uniform $k$-subshift if and only if its
associated \wkautomat on $(\rho \sub {\bf N}(X);\tau _{e_1},\tau _{e_2},\ldots ,\tau _{e_k})$ is an actual
$k$-automaton.

\Proof
  Recall from Diagram \ref{Diagrama} that
  $$
  \tau _n=\rho \sub{\bf N}\theta _n\rho \sub{\bf N}\inv,
  $$
  for every $n$ in ${\bf N}^k$.  Assuming that $X$ is a uniform $k$-subshift, we have that $\rho \sub{\bf N}\inv$ is
uniformly continuous.  Leaving for the reader the easy task of proving that $\theta _n$ is also uniformly continuous, we
deduce that $\tau _n$ is likewise uniformly continuous.  Since we have already seen that $\tau _n$ commutes with the
shift, we deduce from \ref{Hedlund} that $\tau _n$ is a cellular automaton, and in particular so are the $\tau _{e_i}$.
This completes the proof of the ``only if'' part of the statement.

Writing $S_i$ for $\tau _{e_i}$, assume now that each $S_i$ is a cellular automaton, hence a uniformly continuous map.
Setting $Y=\rho \sub{\bf N}(X)$, recall from the proof of \ref{ConstructKSubshift} that the inverse of $\rho \sub{\bf
N}$ on $Y$ is the map $\Psi $ given by
  $$
  \Psi (y)\calcat n = \alpha _n(y)\calcat 0, \for y \in Y, \for n \in {\bf N}^k,
  $$
  where
  $$
  \alpha _n=S_1^{n_1}S_2^{n_2}\cdots S_k^{n_k},
  $$
  for each $n=(n_1,n_2,\ldots ,n_k)\in {\bf N}^k$.

Since each $S_i$, is uniformly continuous, the same is true for $\alpha _n$, so it follows from \ref{UnifCont} that
  $\Psi (y)\calcat n$ depends only on finitely many coordinates of $y$.

Given $\varepsilon >0$, choose $p\in {\bf N}$ such that $2^{-p}<\varepsilon $.  Based on the conclusion of the above
paragraph, let $q$ be a positive integer such that for every $n\in B_p$, one has that $\Psi (y)\calcat n$ depends only
on
  $$
  (y(0),y(1),\ldots y(q)).
  $$
  Setting $\delta =2^{-q}$, assume that $y,y'\in Y$ are such that $d(y,y')\leq \delta $.  Then $y$ and $y'$ agree on
$\{0, 1, \ldots , q\}$, whence
  $$
  \Psi (y)\calcat n = \Psi (y')\calcat n, \for n\in B_p,
  $$
  and we then conclude that
  $$
  d\big( \Psi (y), \Psi (y')\big)\leq 2^{-p}<\varepsilon .
  $$

  This shows that $\Psi $, and hence also $\rho \sub{\bf N}\inv$ is uniformly continuous, which in turn says that $X$ is
a uniform $k$-subshift.  \endProof

Given a $k$-graph $\Lambda $, recall from the paragraph immediately after \ref{DefRankShift}, that the associated ${\bf
N}^k$-subshift $X_\Lambda $ built from $\Lambda $ according to \ref{XLambda} is an example of a $k$-subshift.  Moreover,
as seen in Section \ref{Factoriza}, the \wkautomat on
  $$
  \big(\rho \sub{\bf N}(X_\Lambda );\ S_1,\ldots ,S_k\big)
  $$
  associated to $X_\Lambda $ via \ref{SameThing} is such that the $S_i$ are cellular automata with anticipation 1, and
hence we in fact have a $k$-automaton, according to Definition \ref{DefineGoodStuff.b}.  It then follows from
\ref{CharacStrongAutomata} that $X_\Lambda $ is a uniform $k$-subshift.

It would therefore be interesting to determine conditions on a given $k$-automaton
  $
  (Y;\ S_1,\ldots ,S_k),
  $
  which would imply that it comes from a $k$-graph in the sense that the associated action of ${\bf N}^k$ on $Y$ is
conjugate to one coming from the $k$-automaton arising from a $k$-graph.  By \ref{HighCore}, requiring $Y$ to be a
finite type subshift will likely be among these conditions.

\section Renault-Deaconu groupoids

A submonoid $P$ of a discrete group $G$ satisfying $P\inv P \subseteq P P\inv $ is called an \emph {Ore monoid}.  Given
a right action $\alpha $ of $P$ on a locally compact, Hausdorff, topological space $X$ by means of local homeomorphisms,
as in \cite [Section 2]{ER} (see also \cite {ES}), one may build the corresponding \emph {Renault-Deaconu groupoid}
  $$
  X\rtimes _\alpha P =
  \big \{(x,g,y)\in X\times G\times X : \exists \, n,m\in P,\ \alpha _n(x)=\alpha _m(y),\ g = n m\inv \big \},
  $$
  also called the \emph {transformation}, or \emph {semidirect product} groupoid.  By \cite [3.2]{ER}, one has that
$X\rtimes _\alpha P$ is an \'etale groupoid with the topology generated by the sets of the form
  $$
  U(n,m,A,B) =
    \big \{(x, n m\inv ,y): x\in A,\ y\in B,\ \alpha _n(x)=\alpha _m(y) \big \},
  $$
  where $A,B\subseteq X$ are open subsets, and $n,m\in P$.

\bigskip If $X$ is a given subshift of rank $k$, we of course have a natural action of the Ore monoid ${\bf N}^k$ on
$X$, but this might not put us in the situation of the paragraph above since $k$-subshifts are not necessarily locally
compact, and neither is the action of ${\bf N}^k$ by local homeomorphisms.  In fact, even when $k=1$ and the alphabet is
finite (in which case at least $X$ is compact), namely in the case of a classical subshift, the shift map itself may
fail to be a local homeomorphisms; even worse, it may fail to be an open mapping.  By \cite [2.5]{DE}, only subshifts of
finite type are open.

The $k$-subshift $X_\Lambda $ arising from a row-finite $k$-graph $\Lambda $ fortunately does not suffer from such
tribulations: the path space $\Lambda ^\infty $ is locally compact by \cite [Lemma 2.6]{KP}, and each $\sigma ^p$ is a
local homeomorphisms by \cite [Remarks 2.5]{KP}.  In view of the commutative diagram

\beginmypix x from -400 to 400 y from -80 to 300 xscale 0.012truecm yscale -0.012truecm \def\a{\kern 3pt \Lambda ^\infty
\kern-3pt} \map $\a$ -200 000 $\a$ -200 200 \lbl <-10pt,0pt> {$\sigma ^p$} \map $X_\Lambda $ 000 000 $X_\Lambda $ 000
200 \lbl <-10pt,0pt> {$\theta _p$} \map $\Markov$ 200 000 $\Markov$ 200 200 \lbl <10pt,0pt> {$\tau _p$} \flexa -190 000
000 000 \lbl <0pt,10pt> {$\Xi $} \flexa 000 000 200 000 \lbl <0pt,10pt> {$\rho \sub {\bf N}$} \flexa -190 200 000 200
\lbl <0pt,-10pt> {$\Xi $} \flexa 000 200 200 200 \lbl <0pt,-10pt> {$\rho \sub {\bf N}$} \endmypix

\noindent where we write $\Markov $ for $\rho \sub {\bf N}(X_\Lambda )$, in which all vertical arrows are
homeomorphisms, we then have that $X_\Lambda $ and $\Markov $ are locally-compact spaces, and both $\theta _p$ and $\tau
_p$ are local homeomorphisms. We may therefore form the Renault-Deaconu groupoids relative to the actions $\sigma $,
$\theta $ and $\tau $, obtaining the following three evidently isomorphic groupoids:
  $$
  \Lambda ^\infty \rtimes _\sigma \N , \qquad
  X_\Lambda \rtimes _\theta {\bf N}^k, \and
  \Markov \rtimes _\tau {\bf N}^k.
  $$

The first one above has already explicitly appeared in \cite [Definition 2.7]{KP}, where it was called the \emph {path
groupoid} of $\Lambda $, while playing a prominent role given that its groupoid C*-algebra is isomorphic \cite
[Corollary 3.5]{KP} to the higher rank graph C*-algebra $C^*(\Lambda )$.

Evidently we now see that $C^*(\Lambda )$ may also be modeled by the \wkautomat on $\Markov $.  So let us formally state
this as one of our main conclusions.

\state Theorem
  Given a row-finite $k$-graph $\Lambda $, let $\Sigma $ be the alphabet consisting of all morphisms $\lambda $ with
$d(\lambda )=(1,1,\ldots , 1)$, and let $A=\{A_{\lambda \mu }\}_{\lambda ,\mu \in \Sigma }$ be the 0-1 matrix such that
  $
  A_{\lambda \mu }= 1,
  $
  if and only if $s(\lambda )=r(\mu )$.  Then there are $k$ pairwise commuting cellular automata $S_1,S_2,\ldots ,S_k$,
whose product coincide with Markov's shift on the Markov space $X_A$.  Each $S_i$ is moreover a local homeomorphism and,
denoting by $\tau $ the action of\/ ${\bf N}^k$ on $X_A$ obtained by iterating the $S_i$, one has that the semidirect
product groupoid $X_A\rtimes _\tau {\bf N}^k$ is a model for the higher rank graph C*-algebra in the sense that
  $
  C^*(X_A\rtimes _\tau {\bf N}^k)
  $
  and $C^*(\Lambda )$ are isomorphic.

\references

\Article DE
  M. Dokuchaev and R. Exel;
  Partial actions and subshifts;
  J. Funct. Analysis, 272 (2017), 5038-5106

\Article ER
  R. Exel and J. Renault;
  Semigroups of local homeomorphisms and interaction groups;
  Ergodic Theory Dynam. Systems, 27 (2007), 1737-1771

\Bibitem Kelley
  J. Kelley;
  General topology;
  Springer-Verlag, 1975

\Article KP
  A. Kumjian and D. Pask;
  Higher-rank graph C*-algebras;
  New York J. Math., 6 (2000), 1-20 (electronic)

\Bibitem LindMarcus
  D. Lind and B. Marcus;
  An introduction to symbolic dynamics and coding;
  Cambridge University Press, 1995

\Article ES
  J. Renault and S. Sundar;
  Groupoids associated to Ore semigroup actions;
  J. Operator Theory, 73 (2015), no. 2, 491-514

  \endgroup
  \close